# SIMPLEX AND INTERIOR POINT METHODS FOR SOLVING BUDGETARY ALLOCATION LINEAR PROGRAMMING PROBLEM IN INDUSTRY REVOLUTION 4.0

**ALI KADHIM YAQOOB[1], AHMAD KADRI JUNOH[2], WAN ZUKI AZMAN WAN MUHAMAD[3], MOHAMMAD FADZLI RAMLI[4], NAJAH GHAZALI[5], , MO'ATH ALLUWAICI[6],**

[1,2,3,4,5,6]Institute of Engineering Mathematics, Faculty of Applied and Human Sciences, University Malaysia Perlis, 02600 Arau, Perlis, Malaysia
Email : [2]kadri@unimap.edu.my, [3]wanzuki@unimap.edu.my, [4]mfadzli@unimap.edu.my, [5]najahghazali@unimap.edu.my,

**ABSTRACT**

Increasing the complexity of solving budgetary allocation (NP-hardness problem) has led a wide range of methods to minimize the costs. Metaheuristics and Linear Programming (LP) are the most optimisation in this fields. Therefore, this study provides insights and a deep understanding of the applicability of LP models in industry and how to formulate SM and affine IPM for solving real world linear problems. Moreover, it will present a better way to deal with decision making problems through the development and comparison of the SM and affine IPM to solve LP optimization problem to maximize profit. Finally, to other researchers particularly of similar interests who are undertaking further investigation on this topic, this study can be vital as a secondary source of information and guidance.

## 1. INTRODUCTION

Mathematical models are used extensively in almost all areas of decision-making including profit and budgetary planning, resources allocation.

LP is a mathematical constrained method that is used to solve linear problems particularly to maximize profits and minimize costs subjected to constraints on capacities, resources, demands and supplies . Basically, any linear program constrained model composes of four main parts which are; decision variables, the model parameters, evaluation criterion (objective function) and a set of mathematical constraints. The constraints are mathematical inequalities of equalities which are used to define the solutions to an optimization or constrained mathematical model [21].

Over the past decades, mathematical computational and optimization models particularly Linear Programming (LP) have attempted to provide a systematic, quantitative ways to evaluate and select decisions for budgetary planning in industrial organizations [15]. Budgetary planning is an essential attribute in industries in which a decision on allocating the financial resources may influence the outcome of the profit, cost and the performance of the industry [1].

Recently, various methods have been developed to solve LP problems. However, there is always considerable demand to determine the efficient method to solve large scale LP problems particularly in applied science, engineering and economics. Moreover, the exhausted sophisticated codes which are required to solve linear programming models are quite expensive, time consuming and tedious. Thus, alternative methods are needed to solve large scale linear programming problems such as SM and IPM.

LP problems can be solved through the use of one of two algorithms: either Simplex Method (SM), or Interior Point Methods (IPM)s.

SM was first introduced by George Dantzig in the late 40s to solve large scale LP

problems. It is the very first algorithm for solving LP problems which has been extensively employed for solving LP optimization problems over the past decades in various fields[16]. In order to employ the SM, the LP problem must be represented in a standard form and that the constraints (inequalities) have been converted to equalities using the slack variables as it is not possible to perform SM calculations upon inequality. The standard represents the baseline format for all LP problems before solving for the optimal solution and has three requirements; it must be a maximization problem, its constraints must be in a less than or equal to inequality and all its variables must be non-negative. The slack variables in maximizing problems represent any unused capacity in the constraint and its value can take from zero to the maximum of that constraints. Each constraint has its own separate slack variable. The inequalities in the minimization problems are converted into equalities by subtracting one surplus variable.[17]

IPMs are mathematical techniques specialized in solving large scale LP problems [2]. In mathematical programming, IPMs have been the largest and most dramatic area of research in optimization since the development of the simplex method.[9]

The affine approach is one of the easiest ways to solve the linear system, whether this system is large or small. Over the past decades, affine scaling IPMs have proved reliable in practice in solving large or small LP problems. The affine scaling algorithms consist of three key parts, namely, starting with an interior dual feasible solution, moving towards a better interior solution, and stopping at an optimal dual solution. Moreover, affine scaling works in two phases, the first concerns on finding a feasible point which serves to start optimizing, while the second phase involves conducting the actual optimization while staying strictly to the interior point toward the optimum. [10][12][14][18][19][20][21].

### 1.1 Linear programming (LP)

In today's competitive markets, companies urge to produce high quality products at the lowest possible cost while achieving maximum profit and meeting customer demands. In fact, a company's endurance in a competitive market strictly relies on its ability to manufacture high quality products at lower costs. Based on this, [3] emphasized that, organizations in the world are challenged by shortages of production inputs and low capacity utilization that can consequently lead to low production outputs and costs. [4] addressed the effectiveness of LP as a mathematical optimization tool for determining an optimal solution from a set of an alternative solutions with respect to an evaluation criterion (objective function) and linear mathematical constraints (inequalities or equalities).

LP was introduced as an efficient technique to solve linear problems by Dantzig in 1947 and since then it has been implemented for various practical and industrial [5]. Dantzig defined LP as an extremely revolutionary tool that able to formulate real-world problems in detailed mathematical modeling , to state goals and to select the best decision when facing with a practical situation of great complexity (Dantzig , 1949).

In industries, there are different optimization problems such as to maximize or minimize economic function or a ratio of physical function including cost/volume, cost/time, profit/cost and cost/benefit. Thus, there is a need for developing optimization model to solve these types of problems. Such model is known as the LP model which is expressed [6] as follows:

$$Minimize(orMaximize)\ Z = \sum_{j=1}^{n} c_j x_j$$

Subject to

$$\sum_{j=1}^{n} a_{ij} x_j \{\le, =, \ge\} b_j \; ; i = 1,2,..,m$$

$$x \ge 0; b \ge 0$$

Where, one of the signs $(\le, =, \ge)$ holds for specific constraint and differ from the sign of the other constraints. Here $c_j(j = 1,2,..,n)$ are the coefficients of cost (or profit), $x_j(j = 1,2,..,n)$ are the variables decision. The set of feasible solution to the LP problem (LP) is: $S = \{(x_1, x_2, \dots, x_n)^T : (x_1, x_2, \dots, x_n)^T \in R^2$. In this regard, the $(x_1, x_2, \dots, x_n)^T$ (the set $S$ ) is called the constraints set, feasible set or feasible region of LP model.

### 1.2 Simplex Method (SM)

SM was first introduced by George Dantzig in the late 40s to solve large scale LP problems. It is the very first algorithm for solving

LP problems which has been extensively employed for solving LP optimization problems over the past decades in various fields [7]. It provides an optimal solution for the linear programming problems by examining the vertices of a feasible region through the movement from one extreme point to an adjacent extreme point with lowest possible cost. Basically, simplex represent a mathematical line segment which connect two or more points. For instance, a three-dimensional simplex is a four-sided pyramid which has four corners. SM is an approach to solve LP models using tableaus, variables and pivot variables as means to find the optimal solution of a constrained linear problem [8].

In order to employ the SM, the LP problem must be represented in a standard form and that the constraints (inequalities) have been converted to equalities using the slack variables as it is not possible to perform SM calculations upon inequality. The standard represents the baseline format for all LP problems before solving for the optimal solution and has three requirements; it must be a maximization problem, its constraints must be in a less than or equal to inequality and all its variables must be non-negative. The slack variables in maximizing problems represent any unused capacity in the constraint and its value can take from zero to the maximum of that constraints.

The inequalities in the minimization problems are converted into equalities by subtracting one surplus variable. In fact, standard form is vital as it allows creating an ideal starting point for solving the Simplex method as efficiently as possible as well as other methods for solving LP problems. For instance, Balogun, Jolayemi, Akingbade and Muazu (2012) used SM to solve LP problem to maximize profit from the productions of soft drinks in Bottling Company (Ilorin plant) [17].

### 1.3 Interior Point Method (IPM)

IPMs are mathematical techniques specialized in solving large scale LP problems [2]. In mathematical programming, IPMs have been the largest and most dramatic area of research in optimization since the development of the simplex method [9]. Indeed, IPMs have permanently changed the landscape of mathematical programming theory, practice and computation and overcame the limitation of SM [10].

Generally, IPM solves a LP problems by generating a sequence of points which are inside of the feasible region starting from an initial (strictly) interior point and moving toward the optimum solution while staying closer to the central path away from the boundary regions. The central path of IPMs is basically a smooth curve in the interior feasible region which serves as a guide to the set of optimal solutions. Although there are various types of IPMs, they fall under three categories which are; affine-scaling methods, potential reduction methods and central path methods [11][12]. However, this literature is limited to the affine scale IPM for solving LP problems.

The basic affine scaling algorithm was first developed for LP by a Soviet mathematician named Dikin in 1967 and followed with its convergence in 1974 (Dikin, 1967). Vanderbei et al (1986) replaced the projective transofmation of the Karmarkar by a simpler one known as the affine projective transformation [13]. They proposed the primal affine scaling algorithm to solve the primal linear program problems in standard form and established convergence proofs of the algorithm. However, after a few years, it was realized that the "new" affine scaling algorithms were in fact reinventions of the decades-old results of Dikin. In fact only Barnes and Vanderbei et al (1986) managed to produce an analysis of affine scaling's convergence properties. On the other hand, the dual affine scaling algorithm was formulated by Adler et al (1989) for solving linear programs in inequality form [14]. Though polynomial-time complexity has not been proved yet for this algorithm, global convergence using so-called long steps was proved by Tsuchiya and Muramatsu (1995). This algorithm is often called the primal (or dual) affine-scaling algorithm because the algorithm is based on the primal (or dual) problem only. In the 90s authors have integrated the concepts of dual and primal affine methods and proposed what so called primal-dual affine scaling interior point method for solving LP for example Monteiro, Adler and Resende (1990) and Jansen, Roos, and Terlaky (1996).

The affine scaling algorithms consist of three key parts, namely, starting with an interior dual feasible solution, moving towards a better interior solution, and stopping at an optimal dual solution. Moreover, affine scaling works in two phases, the first concerns on finding a feasible point which serves to start optimizing, while the second phase involves conducting the actual optimization while staying strictly to the interior point toward the optimum.

## 2. METHODOLOGY

The methodology is the most important part of the research and in the eyes of the researcher and summarizes how to choose and determine the method used by and hopes to get through the desired results, as well as the steps on which to reach the solution required. In this paper we use a linear program ready and we apply two methods to see what is the best way between them as well as to find the best solution for this linear program to maximize the profit in it.

### 2.1 Linear programming (LP)

LP is described as a mathematical model used for achieving the maximum or minimum linear function value under specific settings and constraints. In other words, LP is a mathematical method for determining a way to achieve the best outcome (such as maximum profit or minimum cost) in a given mathematical model for some list of requirements represented as linear relationships. Moreover, LP is a technique for the optimization of a linear objective function, subject to linear equality or linear inequality constraints to provide allocation optimality of resources scares as well as maximizing the profit. In this study, the LP is employed to maximize profit of LANA Company based on some constraints and objective function. Therefore, the solution to the LP problem is the set of variables which gives an optimal value to the objective function and at the same time does not violate the constraints imposed upon the variables. Linear programs are problems that can be expressed in canonical form:

Maximize $C^T X$
Subject to $AX \leq b$
where $X \geq 0$

where $X$ represents the vector of variables (to be determined), $C$ and $b$ are vectors of (known) coefficients and $A$ is a (known) matrix of coefficients. The expression to be maximized or minimized is called the objective function ($C^T X$in this case). The equations $AX \leq b$ are the constraints which specify a convex polytope over which the objective function is to be optimized. (In this context, two vectors are comparable when every entry in one is less-than or equal-to the corresponding entry in the other. Otherwise, they are incomparable.).

In this stage, you must have data from the manufacturer or company to maximize its profit, in order to create a linear program by converting the data into different linear equations of equal value. In this paper, the data was collected in advance by researcher as well as the use of his linear program, because this paper shows how to use two methods to solve the linear system, so anyone can learn how to work and convert the data into a linear system through.

**Linear programming:**

$$Max\ Z = 8.073Қ_1 + 6.398Қ_2 + 3.9965\ Қ_3 + 5.943Қ_4 + 5.52175Қ_5 + 7.1955Қ_6$$

**Constraints**

- $Қ_1 + Қ_2 + Қ_3 + Қ_4 + Қ_5 + Қ_6 \geq 74,500$
- $Қ_1 + Қ_2 + Қ_3 + Қ_4 + Қ_5 + Қ_6 \leq 130,000$
- $29.601Қ_1 + 19.194Қ_2 + 21.5811Қ_3 + 22.923Қ_4 + 21.2375Қ_5 + 19.188Қ_6 \geq 1{,}823{,}806.45$
- $8.073Қ_1 + 6.398Қ_2 + 3.9965Қ_3 + 5.943Қ_4 + 5.52175Қ_5 + 7.1955Қ_6 \geq 467{,}663.125$
- $8.073Қ_1 + 6.398Қ_2 + 3.9965Қ_3 + 5.943Қ_4 + 5.52175Қ_5 + 7.1955Қ_6 \leq 765{,}056.25$
- $0.5\ Қ_1 + Қ_2 + 0.5\ Қ_3 + 0.25Қ_4 \leq 50,000$
- $0.25\ Қ_1 + 0.25\ Қ_3 + 0.25\ Қ_4 + 0.5\ Қ_5 \leq 40,000$
- $0.25\ Қ_1 + 0.25\ Қ_3 + 0.5\ Қ_4 + 0.5\ Қ_5 + Қ_6 \leq 40,000$
- $Қ_1 \geq 11,000$
- $Қ_2 \geq 2,200$
- $Қ_3 \geq 8,800$
- $Қ_4 \geq 2,200$
- $Қ_5 \geq 4,400$
- $Қ_6 \geq 2,200$
- $Қ_6 \leq 6{,}500$

### 2.2 Simplex Method (SM)

SM is a basic algorithm for solving LP problems by moving from vertices to vertices along the polytope's edges define the mathematical constraints. The movements of SM algorithms continuous with successive decrease of the values of the evaluation criterion (objective function) until an optimal solution is reached. In this regard, each solution represents an extreme point of the feasible region along the polytope. SM gets to the optimal solution while satisfying the linear problem constraints and optimizes the objective function. During the movement of the SM from one vertex to another (each iteration), the solution in every point is checked and if the optimal solution is not reached the process is repeated using the data obtained from previous iteration. In this regard, the procedures of SM for solving the LP problem are:

i. Converting equations to standard formula by adding variables.
ii. Adding variables to (Z) with multiplying (R) in (-M).
iii. Configure the table as shown in Table 6
iv. Finding the pivot column by estimating the value of the (m) and replacing it with (m) then taking the absolute maximum value negative.
v. Finding the pivot row by dividing the outputs of the equations on the elements of the pivot column and taking the smallest positive value.
vi. Selecting the pivot element which is from the intersection of the pivot row with the pivot column.
vii. Finding the equation of the axial row (A) by dividing the elements of the pivot row on the pivot element.
viii. Find the rest of the rows by applying the following equation:
(the old row - the pivot element in the same row * A).
ix. Putting the new rows in a new table with notes change the pivot row to the (A).
x. Repeat the steps from (4) to (10) until the results are obtained and the values of the function (Z) are zero or positive.

**Linear programming**:

$$Max\ Z = 8.073К_1 + 6.398К_2 + 3.9965\ К_3 + 5.943К_4 + 5.52175К_5 + 7.1955К_6 - 0S1 - 0S2 - 0S3 - 0S4 - 0S5 - 0S6 - 0S7 - 0S8 - 0S9 - 0S10 - 0S11 - 0S12 - 0S13 - 0S14 - 0S15 - MR1 - MR2 - MR3 - MR4 - MR5 - MR6 - MR7 - MR8 - MR9$$

**Constraints**

- $К_1 + К_2 + К_3 + К_4 + К_5 + К_6 - S1 + R1 = 74,500$
- $К_1 + К_2 + К_3 + К_4 + К_5 + К_6 + S2 = 130,000$
- $29.601К_1 + 19.194К_2 + 21.5811К_3 + 22.923К_4 + 21.2375К_5 + 19.188К_6 - S3 + R2 = 1{,}823{,}806.45$
- $8.073К_1 + 6.398К_2 + 3.9965К_3 + 5.943К_4 + 5.52175К_5 + 7.1955К_6 - S4 + R3 = 467{,}663.125$
- $8.073К_1 + 6.398К_2 + 3.9965К_3 + 5.943К_4 + 5.52175К_5 + 7.1955К_6 + S5 = 765{,}056.25$
- $0.5\ К_1 + К_2 + 0.5\ К_3 + 0.25К_4 + S6 = 50,000$
- $0.25\ К_1 + 0.25\ К_3 + 0.25\ К_4 + 0.5\ К_5 + S7 = 40,000$
- $0.25\ К_1 + 0.25\ К_3 + 0.5\ К_4 + 0.5\ К_5 + К_6 + S8 = 40,000$
- $К_1 - S9 + R4 = 11,000$
- $К_2 - S10 + R5 = 2,200$
- $К_3 - S11 + R6 = 8,800$
- $К_4 - S12 + R7 = 2,200$
- $К_5 - S13 + R8 = 4,400$
- $К_6 - S14 + R9 = 2,200$
- $К_6 + S15 = 6{,}500$

### 2.3 Primal Affine Interior Point Method

The method of affine is one of the easiest ways to solve the linear system, whether this system is large or small. This method is based on a number of clear and specific steps to solve the linear system. This method depends on converting the linear system into matrices with initial values called ($X_0$). The solution is initiated by these values to find other values that are added to the linear system regression to be converted into equations. The values used and the values that were found consist of the first matrix in the solution. The rest of the steps are completed until the values ($X_1$) are obtained. Finally, the stopping condition (access to the solution) is the values that were found after the implementation of all steps are values ($X_1$). If the optimal solution is not reached then, IPMs move on to the second step and repeat the same steps.

**Steps of the Primal Affine Algorithm:**

a. Conversion of the inequations into equations by adding the value of each of

the assumptions so that this value achieves the required equality when the variable is added. Thus, if the inequations are (≥) then the (-) is added and if they (≤) then the (+) is added.

$$Max\ Z = 8.073К_1 + 6.398К_2 + 3.9965\ К_3 + 5.943К_4 + 5.52175К_5 + 7.1955К_6$$

Subject to:

- $К_1 + К_2 + К_3 + К_4 + К_5 + К_6 - S_7 = 74,500$
- $К_1 + К_2 + К_3 + К_4 + К_5 + К_6 + S_8 = 130,000$
- $29.601К_1 + 19.194К_2 + 21.5811К_3 + 22.923К_4 + 21.2375К_5 + 19.188К_6 - S_9 = 1{,}823{,}806.45$
- $8.073К_1 + 6.398К_2 + 3.9965К_3 + 5.943К_4 + 5.52175К_5 + 7.1955К_6 - S_{10} = 467{,}663.125$
- $8.073К_1 + 6.398К_2 + 3.9965К_3 + 5.943К_4 + 5.52175К_5 + 7.1955К_6 + S_{11} = 765{,}056.25$
- $0.5\ К_1 + К_2 + 0.5\ К_3 + 0.25К_4 + S_{12} = 50,000$
- $0.25\ К_1 + 0.25\ К_3 + 0.25\ К_4 + 0.5\ К_5 + S_{13} = 40,000$
- $0.25\ К_1 + 0.25\ К_3 + 0.5\ К_4 + 0.5\ К_5 + К_6 + S_{14} = 40,000$
- $К_1 - S_{15} = 11,000$
- $К_2 - S_{16} = 2,200$
- $К_3 - S_{17} = 8,800$
- $К_4 - S_{18} = 2,200$
- $К_5 - S_{19} = 4,400$
- $К_6 - S_{20} = 2,200$
- $К_6 + S_{21} = 6{,}500$

b. The values obtained from solving LP by using SM are used to obtain variables that have been added, for example: $К_1 + К_2 + К_3 + К_4 + К_5 + К_6 - S_7 = 74,500$

From the SM expected result:

- $40053 + 16758 + 8801 + 2200 + 48971 + 2201 - S_7 = 74,500$
- $118984 - S_7 = 74,500$
- $-S_7 = 74,500 - 118984 = -44484$
- $S_7 = -44484$ By adding ($S_7$) value the following can be obtained:
- $40053 + 16758 + 8801 + 2200 + 48971 + 2201 - 44484 = 74,500$

**c.** The corresponding feasible interior point starting solution is $x^0$ , $X^0 = \begin{bmatrix} K_1\ K_2\ K_3\ K_4\ K_5\ K_6\ S_7\ S_8\ S_9\ S_{10}\ S_{11}\ S_{12} \\ S_{13}\ S_{14}\ S_{15}\ S_{16}\ S_{17}\ S_{18}\ S_{19}\ S_{20}\ S_{21} \end{bmatrix}^T$

d. The first iteration started by setting the scaling matrix D (21×21 ) which is $X^0$ in step 3.

e. Matrix (A) consisting of equations by taking variables coefficients in the step 1

f. Matrix (Â); $\hat{A} = A * D$

g. Matrix (C) consisting of Variable coefficients (Max Z)

h. Matrix (Ĉ) and $\hat{C} = D * C$

i. Matrix $p = I - \hat{A}^t * \left(\hat{A} * A^t\right)^{-1} * \hat{A}$

j. Matrix $p^{\circ}$; $(p^0 = -p * c)$

k. For (ϴ) if $(p^0 = -p * c)$ then (ϴ) is the negative number that have the largest absolute value.

l. Matrix $\dot{X}^1$ , $\dot{X}^1 = e - \frac{\alpha}{\theta} * p^t$ and the $(\alpha)$=0.5
In the affine scaling interior point algorithm the selected constant, $(\alpha)$ is required to be such that 0< $\alpha$ <1 typically, α is set to be between 0.5 and 0.95 for primal affine methods (Dikin, 1967). Dikin converged the α using $\alpha$ = 0.5 and since then it was widely used in literature by various researchers.

m. Matrix $X^1$, $X^1 = D * X$. At this point, the first iteration of the method is completed and the comparison of the matrix $X^1$ with the matrix $X^0$ is drawn. If the results are similar then the stopping criterion is met and this is the optimal solution. Meanwhile, if the results are not identical and there is no approximation between them, then the model proceeds to the next iteration which involves similar steps of the solution except the scaling matrix (D) whose values will change from $X^0$ to $X^1$ and the solution will be repeated until the optimal solution is obtained.

$$D = \begin{bmatrix}
40053 & 0 & 0 & 0 & 0 & 0 & 0 & 0 & 0 & 0 & 0 & 0 & 0 & 0 & 0 & 0 & 0 & 0 & 0 & 0 & 0 \\
0 & 16758 & 0 & . & . & . & . & . & . & . & . & . & . & . & . & . & . & . & . & 0 & 0 \\
0 & 0 & 8801 & 0 & . & . & . & . & . & . & . & . & . & . & . & . & . & . & . & 0 & 0 \\
. & . & 0 & 2200 & 0 & . & . & . & . & . & . & . & . & . & . & . & . & . & . & 0 & 0 \\
. & . & . & 0 & 48971 & 0 & . & . & . & . & . & . & . & . & . & . & . & . & . & 0 & 0 \\
. & . & . & . & 0 & 2200 & 0 & . & . & . & . & . & . & . & . & . & . & . & . & 0 & 0 \\
. & . & . & . & . & 0 & -44484 & 0 & . & . & . & . & . & . & . & . & . & . & . & 0 & 0 \\
. & . & . & . & . & . & 0 & 11016 & 0 & . & . & . & . & . & . & . & . & . & . & 0 & 0 \\
. & . & . & . & . & . & . & 0 & -2006075.71 & 0 & . & . & . & . & . & . & . & . & . & 0 & 0 \\
. & . & . & . & . & . & . & . & 0 & -301178.339 & 0 & . & . & . & . & . & . & . & . & 0 & 0 \\
. & . & . & . & . & . & . & . & . & 0 & 0.01425 & 0 & . & . & . & . & . & . & . & 0 & 0 \\
. & . & . & . & . & . & . & . & . & . & 0 & 8265 & 0 & . & . & . & . & . & . & 0 & 0 \\
. & . & . & . & . & . & . & . & . & . & . & 0 & 2751 & 0 & . & . & . & . & . & 0 & 0 \\
. & . & . & . & . & . & . & . & . & . & . & . & 0 & 0 & 0 & . & . & . & . & 0 & 0 \\
. & . & . & . & . & . & . & . & . & . & . & . & . & 0 & -29053 & 0 & . & . & . & 0 & 0 \\
. & . & . & . & . & . & . & . & . & . & . & . & . & . & 0 & -14558 & 0 & . & . & 0 & 0 \\
. & . & . & . & . & . & . & . & . & . & . & . & . & . & . & 0 & 1 & 0 & . & 0 & 0 \\
. & . & . & . & . & . & . & . & . & . & . & . & . & . & . & . & 0 & 0 & 0 & 0 & 0 \\
. & . & . & . & . & . & . & . & . & . & . & . & . & . & . & . & . & 0 & -44571 & 0 & 0 \\
0 & 0 & 0 & 0 & 0 & 0 & 0 & 0 & 0 & 0 & 0 & 0 & 0 & 0 & 0 & 0 & 0 & 0 & 0 & 1 & 0 \\
0 & 0 & 0 & 0 & 0 & 0 & 0 & 0 & 0 & 0 & 0 & 0 & 0 & 0 & 0 & 0 & 0 & 0 & 0 & 0 & 4299
\end{bmatrix}$$

## 3. RESULTS AND DISCUSSION

Present and elaborate the results and discussion of the Linear Programming (LP), Simplex Method (SM) and Interior Point Method (IPM) models for solving budgetary planning in industries. To demonstrate the applicability of the SM and IPM for solving LP profit, a case study was presented which is to maximize the profit of LANA company for Food Ltd. Besides that, this section presents a comparison between the simplex method (SM) and the interior point methods (IPM) for solving linear programming (LP). It also illustrates the expected results.

### 3.1. Linear programming

Linear Programming (LP) is as a mathematical constrained method, which can be utilized for achieving the maximum or minimum linear function value under specific settings and constraints. It is a dominant technique for providing allocation optimality of resources scares as well as maximizing the profit. It is employed for optimizing the profit of LANA company for Food Ltd subjected to a particular a set of related variables that are independent in a relationship (linear constraints). Here, the profit to be maximized can be designated as the dependent variable which can be represented by the objective function. Meanwhile, the independent variables are unknown values which are determined by solving the linear problem. In this study, LP model was successfully developed with the purpose of maximizing the profit of LANA company for Food Ltd. Toward this end, an objective function was formulated using six products namely Red Beans ($K_1$), Green Beans ($K_2$), Chick Peas ($K_3$), Hamos ($K_4$), White Beans ($K_5$) and Large Beans ($K_6$). The LP model was subjected to various constraints. Table 1 shows the results of LP optimization program in comparison with the model constraints. The results indicated that for the LP model the optimal global solution is found to be 765,056.25 SAR. Thus, the factory is required to achieve optimal profit of 765,056.25 SAR

*Table 1: Comparison between achieved values for products and the constraints.*

| Variable (product) | Obtained Values | Constraints |
|---|---|---|
| $K_1$ | 78,525 | $X_1 \geq 11,000$ |
| $K_2$ | 2,200 | $X_2 \geq 2,200$ |
| $K_3$ | 15,974 | $X_3 \geq 8,800$ |
| $K_4$ | 2200 | $X_4 \geq 2,200$ |
| $K_5$ | 4,400 | $X_5 \geq 4,400$ |
| $K_6$ | 2,201 | $2,200 \leq X_6 \leq 6,500$ |

### 3.2 Simplex Method

Simplex Method (SM) was developed to further maximize the profit of LANA company for Food Ltd based on the LP for the six products namely Red Beans ($K_1$), Green Beans ($K_2$),

Chick Peas ($K_3$) , Hamos ($K_4$), White Beans ($K_5$) and Large Beans ($K_6$) for 10 months in the year 2014. Here, SM algorithm solve LP problems by moving from vertices to vertices along the polytope's edges defined by the mathematical constraints. The movements of SM algorithms continuous with successive decrease of the values of the evaluation criterion (objective function) until an optimal profit of LANA company for Food Ltd is reached. we used the (WINQSB) and (QM) programing to get the results in (SM) method because the (LP) is big, so it is hard do it by hand but if the (LP) is small can used the hand way by steps in the methodology.

The optimal profit obtained using the Simple Method With values of [765005 by (WinQsb ) , 765056 by (QM)] SAR. Meanwhile, Table 2 depicts the results of the SM. It further compares the obtained constraints with the original constraints. It be clearly observed that the developed simplex method has successfully achieved the constraints and provided the optimal solution.

*Table 2: Comparison between achieved values for products and the constraints using simplex method (winQsb) & (QM) programming.*

| Variable (product) | Obtained Values winQsb | Obtained Values QM | Constraints |
|---|---|---|---|
| $K_1$ | 40053 | 22856 | $X_1 \geq 11,000$ |
| $K_2$ | 16750 | 33619 | $X_2 \geq 2,200$ |
| $K_3$ | 8801 | 8800 | $X_3 \geq 8,800$ |
| $K_4$ | 2200 | 2200 | $X_4 \geq 2,200$ |
| $K_5$ | 48971 | 54579 | $X_5 \geq 4,400$ |
| $K_6$ | 2201 | 2200 | $2,200 \leq X_6 \leq 6,500$ |

### 3.3 Interior Affine Point Method

The affine point method was developed to determine the optimal maximum profit for LANA company for Food Ltd based on the LP for the six products namely Red Beans ($K_1$), Green Beans ($K_2$) , Chick Peas ($K_3$) , Hamos ($K_4$), White Beans ($K_5$ ) and Large Beans ($K_6$ ) for 10 months in the year 2014. Unlike Simplex Method, the IPM method approaches the optimum from the interior of the feasible solution space at the boundary of the feasible region. Thus, IPM solve linear programming by iterating from the interior of the polytope defined by the constraints.

Based on the results of the affine Interior Point Method, the profit obtained is 765289.9244 SAR when we used the SM (WinQsb) result program . Also based on the results of the affine Interior Point Method when we used the SM (QM) result program by same affine Interior Point Method steps , the profit obtained is 765121.8775 . In addition, the resultant constraints of the affine scaling interior point method are shown in Table 3. It be clearly observed that the developed IPM has successfully achieved the constraints and provided the optimal solution.

*Table 3: Comparison between achieved values for products and the constraints using affine- scaling IPM.*

| Variable (product) | Obtained Values 1 | Obtained Values 2 | Constraints |
|---|---|---|---|
| $K_1$ | 40054 | 16938 | $X_1 \geq 11,000$ |
| $K_2$ | 16755 | 36582 | $X_2 \geq 2,200$ |
| $K_3$ | 8805.3 | 8800 | $X_3 \geq 8,800$ |
| $K_4$ | 2200 | 2200 | $X_4 \geq 2,200$ |
| $K_5$ | 49013 | 5981 | $X_5 \geq 4,400$ |
| $K_6$ | 2200.4 | 2200 | $2,200 \leq X_6 \leq 6,500$ |

### 3.4 Comparison between SM and IPM

Interior point algorithm is a polynomial time algorithms. This means that the time required to solve an LP problem of size n would take at most $an^b$ where a and b are two positive numbers. On the other hand, the Simplex algorithm is an exponential time algorithm in solving LP problems. This implies that, in solving an LP problem of size n there exists a positive number such that for any of the Simplex algorithm would find its solution in a time of at most c2^n. For large enough n (with positive a, band c), $c2^n > an^b$ . This means that, in theory, the polynomial time algorithms are superior to exponential algorithms for large LP problems.

Based on the results obtained for both SM and affine IPM, it can be observed that the IPM has produced better results for profit than SM. With values of [765005 by (WinQsb ) , 765056 by (QM)] SAR and is [765289.9244 , 765121.8775] SAR for SM and IPM respectively. Furthermore, the results of the constraints (six products) obtained using affine IPM are better than those obtained using SM as elaborated in Table 4. The (IPM) depends on the (SM) in the solution. It is not possible to start the solution if there are no initial values. These values either are set by the factory or are found by solving the linear system by (SM).

*Table 4: Comparison between achieved values for products and the constraints using of SM and IPM.*

| Variable (product) | Obtained SM(winQsb) | Obtained SM (QM) | Obtained (IPM)1 | Obtained (IPM)2 | Constraints |
|---|---|---|---|---|---|
| $K_1$ | 40053 | 22856 | 40054 | 16938 | $X_1 \geq 11,000$ |
| $K_2$ | 16750 | 33619 | 16755 | 36582 | $X_2 \geq 2,200$ |
| $K_3$ | 8801 | 8800 | 8805.3 | 8800 | $X_3 \geq 8,800$ |
| $K_4$ | 2200 | 2200 | 2200 | 2200 | $X_4 \geq 2,200$ |
| $K_5$ | 48971 | 54579 | 49013 | 59810 | $X_5 \geq 4,400$ |
| $K_6$ | 2201 | 2200 | 2200.4 | 2200 | $2,200 \leq X_6 \leq 6,500$ |

Overall findings indicated that the affine scaling IPM produced optimal profit than the Simplex Method for LANA company for Food Ltd for six products namely Red Beans ($К_1$), Green Beans ($К_2$), Chick Peas ($К_3$), Hamos ($К_4$), White Beans ($К_5$) and Large Beans ($К_6$) for a period of 10 months in the year 2016.

It can be simply observed that the IPM1 and IPM2 produced better results than SM (winQsb) and SM (QM) methods. This is can be clearly observed by the Analysis of Variance (ANOVA) test results in Table 5 and Table 6 and Figure 3.1 and Figure 3.2. The ANOVA tests results are divided into two models: The first ANOVA model comprises of LP, SM (QM) and IPM1. Meanwhile, the second ANOVA model contains LP, SM (winQsb) and IPM2. In both cases the SPSS and Mintab software were used for the analysis and to obtain the results of the ANOVA models. The LP model was chosen as the dependent variables and the SM and IPM are the independent variables as seen in Table 5 and Table 6.

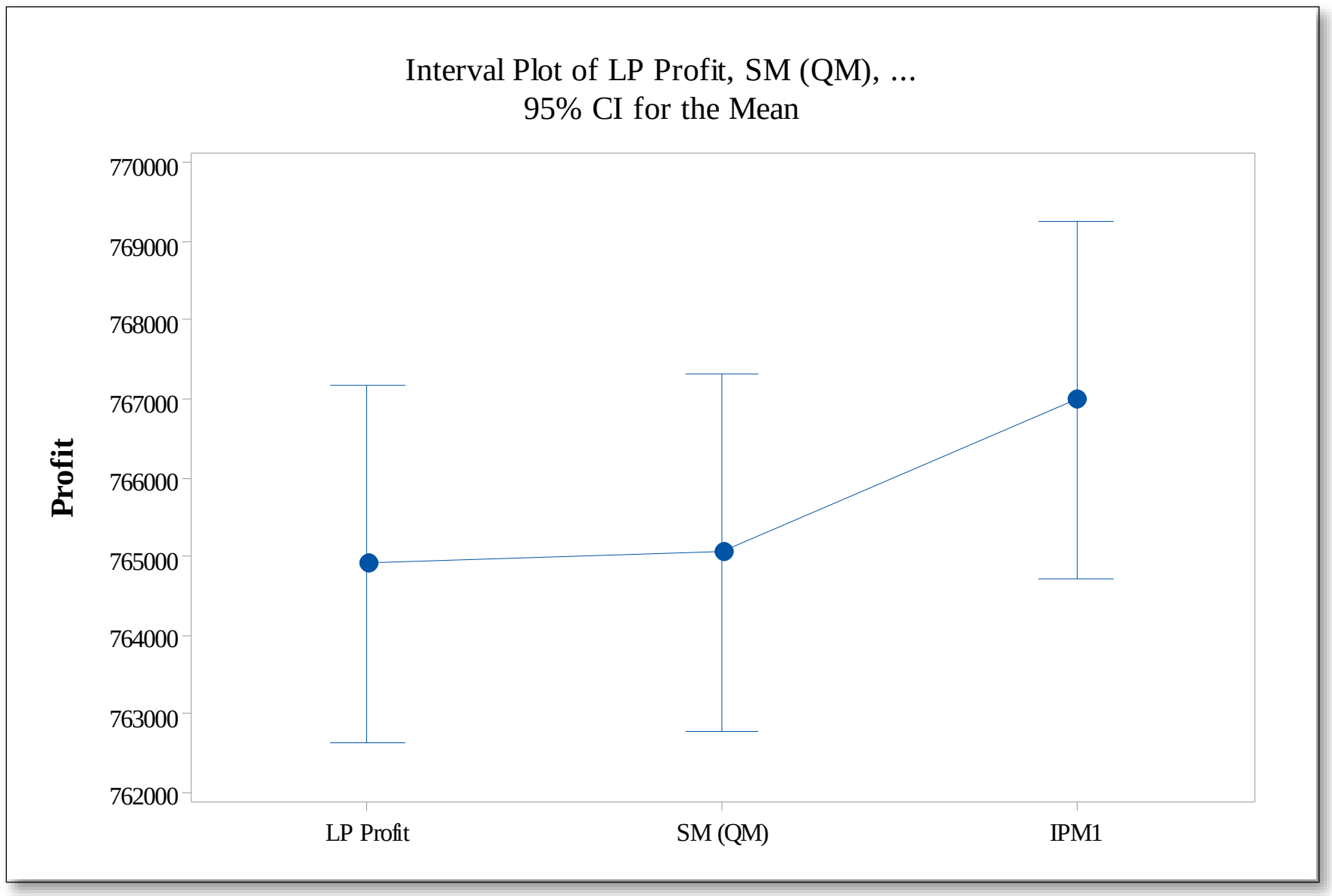

*Figure 3.1: Analysis of variance for the resultant profits of LP, SM (QM) and IPM1 models*

*Table 5: One way ANOVA test results for LP, SM (QS) and IPM1 models*

| | | Sum of Squares | df | Mean Square | F | Sig. |
|---|---|---|---|---|---|---|
| **SM (QS)** | Between Groups | 0.105 | 5 | 0.021 | 0.840 | 0.583 |
| | Within Groups | 0.100 | 4 | 0.025 | | |
| | Total | 0.205 | 9 | | | |
| **IPM1** | Between Groups | 84700537.720 | 5 | 16940107.540 | 0.278 | 0.904 |
| | Within Groups | 243694247.400 | 4 | 60923561.840 | | |
| | Total | 328394785.100 | 9 | | | |

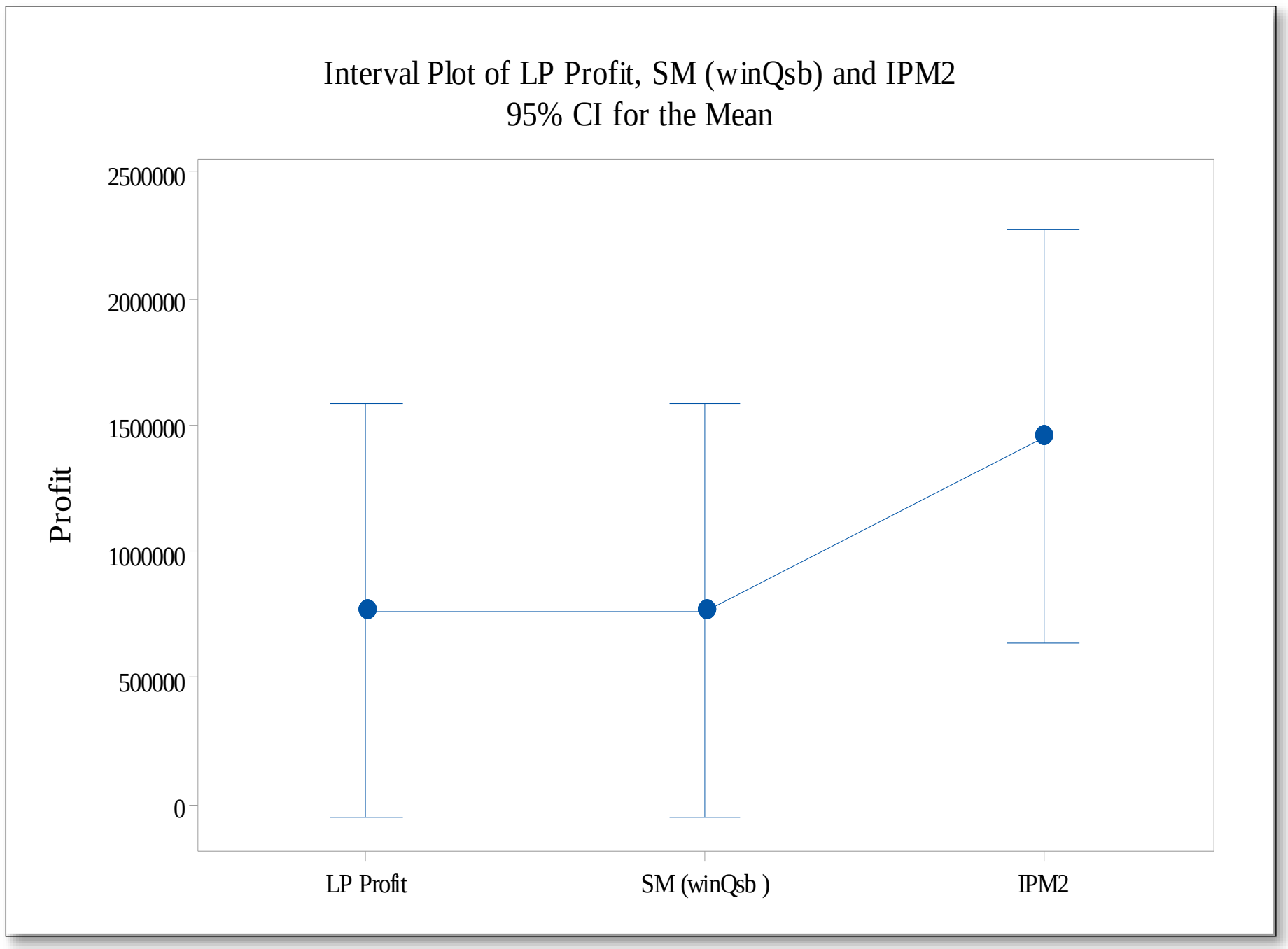

*Figure 3.2: Analysis of variance for the resultant profits of LP, SM (winQsb) and IPM2 models*

*Table 6: One way ANOVA test results for LP, SM (winQsb) and IPM2 models*

| | | Sum of Squares | df | Mean Square | F | Sig. |
|---|---|---|---|---|---|---|
| SM (winQsb) | Between Groups | 0.237 | 5 | 0.047 | 28.480 | 0.003 |
| | Within Groups | 0.007 | 4 | 0.002 | | |
| | Total | 0.244 | 9 | | | |
| IPM2 | Between Groups | 11082823730000.000 | 5 | 2216564745000 | 0.281 | 0.902 |
| | Within Groups | 31577294760000.000 | 4 | 7894323689000 | | |
| | Total | 42660118480000.000 | 9 | | | |

## 4. CONCLUSION

To recap, the first objective of this study was to develop linear programing model to maximize the profit of LANA Company for Food Ltd. The second objective was to develop Simplex Method to maximize the profit of LANA Company for Food Ltd. The third objective was to develop interior point method to further maximize the profit of LANA Company for Food Ltd. Finally, the fourth objective was to compare between the obtained profits with the actual profit of the LANA Company for Food Ltd.

At first the LP mathematical constraint model was developed with the objective function to maximize profit of LANA Company for Food Ltd. The LP model was developed using the LINGO software based on the objective function and subjected to the model constraints. LP has effective capability for solving linear problems particularly profit maximization in industries. The findings showed that the developed LP provided maximum profit when compared to the actual profit.

Apart from this, Simplex Method (SM) was developed to maximize the profit of the aforementioned company based on the LP model. The SM is well known efficient method for solving various linear programming problems. The results obtained from the SM methods were then utilized to develop the Affine-scaling Interior Point Method (IPM) optimization model. The WINQSB and QM software was used as a platform to develop the algorithm for the SM. The optimal profit of LANA company was achieved considering six products namely Red Beans ($K_1$), Green Beans ($K_2$), Chick Peas ($K_3$), Hamos ($K_4$), White Beans ($K_5$) and Large Beans ($K_6$).

Overall conclusive remarks in this study indicated that the affine scaling IPM produced optimal profit than the Simplex Method for LANA company for Food Ltd for six products namely Red Beans ($K_1$), Green Beans ($K_2$), Chick Peas ($K_3$), Hamos ($K_4$), White Beans ($K_5$) and Large Beans ($K_6$) for a period of 10 months in the year 2016. The obtained optimal profits are [765005 by (WinQsb ) , 765056 by (QM)] SAR and [765289.9244 , 765121.8775] SAR for SM and IPM respectively. Finally, the obtained results of both SM and IPM are compared with the actual profit of the industry. The results presented in this research to maximize the profit outperformed the actual profit of the industry. Thus, the budgetary optimization methods demonstrated in this research can serves as a guidance to the industry toward achieving maximum profit while maintaining productivity and high-quality products.

## AKNOWLEDGEMENT

This study was financially supported by the Ministry of Higher Education, Malaysia, under the Fundamental Research Grant Scheme (FRGS/1/2019/TK03/UNIMAP/02/7)